\documentclass{article}
\usepackage{amsmath,amscd}
\usepackage{amssymb}
\usepackage{amsthm}
\newtheorem{theorem}{Theorem}
\newtheorem{question}[theorem]{Question}
\newtheorem{lemma}[theorem]{Lemma}

\newtheorem{definition}[theorem]{Definition}

\newtheorem{remark}[theorem]{Remark}

\newtheorem{desired result}[theorem]{Desired result}

\newtheorem{sublemma}[theorem]{Sublemma}
\newtheorem{notation}[theorem]{Notation}

\newcommand{\begd}{\begin{displaystyle}}

\newcommand{\ga}{\alpha}
\newcommand{\gb}{\beta}

\newcommand{\gD}{\Delta}

\newcommand{\mbq}{\mathbb{Q}}

\newcommand{\mbz}{\mathbb{Z}}

\newcommand{\gal}{\textrm{Gal}\,}
\newcommand{\lcm}{\textrm{lcm}\,}
\newcommand{\ol}[1]{\overline{#1}}
\newcommand{\mc}[1]{\mathcal{#1}}

\title{A bound for the ``torsion conductor'' of a non-CM elliptic curve}
\author{Nathan Jones}
\date{}

\begin{document}
\maketitle 

\begin{abstract}
Given a non-CM elliptic curve $E$ over $\mbq$, define the ``torsion conductor'' $m_E$ to be the smallest positive integer
so that the Galois representation on the torsion of $E$ has image $\pi^{-1}(\gal(\mbq(E[m_E])/\mbq)$, where
$\pi$ denotes the natural projection $GL_2(\hat{\mbz}) \rightarrow GL_2(\mbz/m_E\mbz)$.  We show
that, uniformly for semi-stable non-CM elliptic curves $E$ over $\mbq$, one has ${m_E \ll \left( \prod_{p \mid\gD_E} p\right)^5}$.
\end{abstract}


\section{Introduction}

Let $E$ be an elliptic curve defined over a number field $K$ and let
\[
\varphi_{E} : \gal(\ol{K}/K) \rightarrow GL_2(\hat{\mbz})
\]
be the continuous group homomorphism defined by letting $\gal(\ol{K}/K)$ operate on the
torsion points of $E$ and by choosing an isomorphism $\text{Aut}(E_{\text{tors}}) \simeq GL_2(\hat{\mbz})$.
We will refer to $\varphi_E$ as the \textbf{torsion representation of $E$}.
A celebrated theorem of Serre \cite{serre} shows that, if $E$ has no complex multiplication, then
the index of the image of $\varphi_E$ is finite:
\[
[GL_2(\hat{\mbz}) : \varphi_E(\gal(\ol{K}/K))] < \infty.
\]
This is equivalent to the statement that there exists an integer $m \geq 1$ with the property that
\begin{equation} \label{inverseimage}
\varphi_E(\gal(\ol{K}/K)) = \pi^{-1}(\gal(K(E[m])/K)),
\end{equation}
where $K(E[m])$ denotes the $m$-th division field of $E$, obtained by adjoining to $K$ the $x$ and $y$
coordinates of the $m$-torsion points of a Weierstrass model of $E$, and
\[
\pi : GL_2(\hat{\mbz}) \rightarrow GL_2(\mbz/m\mbz)
\]
denotes the projection.
\begin{definition}
We define the \textbf{torsion conductor} $m_E$ of a non-CM elliptic curve $E$ over $K$ to be the smallest positive
integer $m$ so that \eqref{inverseimage} holds.
\end{definition}
Serre \cite[p. $299$]{serre} has asked the following important question about the image of $\varphi_E$.
\begin{question} \label{serresquestion} Given a number field $K$, is there a constant $C_K$,
such that, for any non-CM elliptic curve $E$ over $K$ and any rational prime number $p \geq C_K$, one has
\[
\gal(K(E[p])/K) \simeq GL_2(\mbz/p\mbz)?
\]
\end{question}
Even in the case of $K = \mbq$ this question remains unanswered.  Mazur \cite[Theorem $4$, p. 131]{mazur} has shown that, 
\begin{equation} \label{mazursresult}
\text{$E$ is semi-stable } \; \Longrightarrow \; \forall \, p \geq 11, \;  \gal(\mbq(E[p])/\mbq) = GL_2(\mbz/p\mbz)
\end{equation}
His work also shows that,
if $p > 19$, 
$p \notin \{ 37, 43, 67, 163 \}$, and
\begin{equation} \label{pexceptional}
\gal(\mbq(E[p])/\mbq) \subsetneq GL_2(\mbz/p\mbz),
\end{equation}
then $\gal(\mbq(E[p])/\mbq)$ is contained in the normalizer of a Cartan subgroup of $GL_2(\mbz/p\mbz)$.  The work of Parent \cite{parent} represents further progress towards resolution of the split Cartan case, while the work of Chen \cite{chen} shows that in the non-split case, new ideas are needed.
Other authors have bounded the largest prime $p$ satisfying \eqref{pexceptional} in terms of invariants of the elliptic curve (\cite{serre3}, \cite{kraus}, \cite{cojocarukani}, and \cite{mw}).

In some applications it is useful to 
have effective control over the variation of $m_E$ with $E$.  For example, in \cite{jonesltn},
such control becomes necessary to compute averages of various constants attached to elliptic curves.
In this note we prove the following theorem.

\begin{theorem} \label{main}
Let $\gD_E$ denote the minimal discriminant of an elliptic curve $E$ over $\mbq$.  Then, uniformly for semi-stable non-CM elliptic curves $E$ over $\mbq$, one has
\[
m_E \ll \left(\prod_{p \text{ prime}, \; p \mid \gD_E} p \right)^5.
\]
If Question \ref{serresquestion} has an affirmative answer when $K = \mbq$, then the above bound holds uniformly for all elliptic curves $E$ over $\mbq$.
\end{theorem}

The proof of Theorem \ref{main} uses elementary Galois theory to reduce the question to working ``vertically over exceptional primes'',
or in other words, to the analogous question of the Galois representation on the Tate module
\[
\gal(\ol{\mbq}/\mbq) \rightarrow GL_2(\mbz_p),
\]
where $p$ satisfies \eqref{pexceptional}.  Such a study has
been carried out in the recent work of Arai \cite{arai}.  The main ideas are present in \cite{serre2} and \cite{langtrotter}.
\begin{remark}
The torsion conductor $m_E$ should not be confused with the number
\[
A(E) := 2\cdot 3 \cdot 5 \cdot \prod_{{\begin{substack} {p \text{ prime } \\ \gal(\mbq(E[p])/\mbq) \subsetneq GL_2(\mbz/p\mbz) } \end{substack}}} p,
\]
discussed in \cite{cojocarukani}, which has the useful property that, for any integer $n$,
\[
\gcd(n,A(E)) = 1 \; \Longrightarrow \; \gal(\mbq(E[n])/\mbq) = GL_2(\mbz/n\mbz).
\]
This condition is weaker than \eqref{inverseimage}.  For example, if $E$ is the curve $y^2 + y = x^3 - x$, then $A(E) = 30$ and $m_E = 74$.  More generally, when $E$ is a Serre curve (for a definition, see \cite[pp. $310$--$311$]{serre} or \cite[Section $3$]{jonesltn}), one has $A(E) = 30$, whereas $m_E$ is greater than or equal to the square-free part of $|\gD_E|$\footnote{By the square-free part $|\gD_E|$, we mean the unique square-free number $n$ such that $|\gD_E|/n$ is a square.}.
\end{remark}
\begin{notation}
For a fixed elliptic curve $E$ over $\mbq$ and for any positive integer $n$ we will denote
\[
L_{n} := \mbq(E[n]), \quad G(n) := \gal(L_{n}/\mbq),
\]
and we will regard $G(n)$ as a subgroup of $GL_2(\mbz/n\mbz)$.  Also, we will
overwork the symbol $\pi$, using it to denote any one of the canonical projections
\[
\begin{split}
&\pi : GL_2(\hat{\mbz}) \rightarrow GL_2(\mbz/n\mbz), \quad \pi : GL_2(\mbz_p) \rightarrow GL_2(\mbz/p^n\mbz), \\
\text{or} \quad &\pi : GL_2(\mbz/n\mbz) \rightarrow GL_2(\mbz/d\mbz) \quad (d \text{ dividing } n),
\end{split}
\]
or the restrictions of any of these projections to closed subgroups, for example
\[
\pi : \varphi_E(\gal(\ol{\mbq}/\mbq)) \rightarrow G(M) \quad \text{ or } \quad \pi : G(n) \rightarrow G(d) \quad (d \text{ dividing } n).
\]
In ambiguous instances, we will denote alternatively
\[
\pi_{n,d} : GL_2(\mbz/n\mbz) \rightarrow GL_2(\mbz/d\mbz).
\]
\end{notation}
We hope that these abbreviations will minimize cumbersome notation and not cause any confusion.
We will say that an integer $M$ \emph{divides} $N^{\infty}$ if whenever a prime $p$ divides $M$, $p$ also divides $N$.  Throughout, the letters $p$ and $\ell$ will always denote prime numbers.

\section{Acknowledgments}

I would like to thank C. David and A. C. Cojocaru for stimulating discussions and for comments on an earlier version.

\section{Proof of Theorem \ref{main}}

Let $E$ be a fixed non-CM elliptic curve over $\mbq$ and denote by
\[
\varphi_{E,p} : \gal(\ol{\mbq}/\mbq) \rightarrow GL_2(\mbz_p) \simeq \text{Aut}(\lim_{\leftarrow} E[p^n])
\]
the Galois representation on the Tate module of $E$ at $p$.  The following is a re-statement of \cite[Theorem $1.2$]{arai}.
\begin{theorem} \label{araisthm}
Let $K$ be a number field and let $p$ be a prime number.  There exists an exponent $n_K(p)$ so that, for each non-CM elliptic curve $E$ over $K$
one has
\begin{equation*} \label{stable}
\varphi_{E,p}(\gal(\ol{K}/K)) = \pi^{-1}(\gal(K(E[p^{n_K(p)}])/K)).
\end{equation*}
\end{theorem}
If $n_K(p) = 0$, this is interpreted to mean that $\varphi_{E,p}$ is surjective.  In fact, for $p > 3$ one has
\begin{equation} \label{tatemodulesurjection}
G(p) = GL_2(\mbz/p\mbz) \; \Longrightarrow \; n_\mbq(p) = 0.
\end{equation}
This is proved by applying \cite[Lemma $3$, p. IV-23]{serre2} with $X$ equal to the commutator subgroup of $\varphi_{E,p}(\gal(\ol{\mbq}/\mbq))$, together with
the fact that because of the Weil pairing, the determinant map
\[
\det : \gal(L_{p^\infty}/\mbq) \twoheadrightarrow (\mbz_p)^*
\]
is surjective.  We define
\[
S := \{ 2, 3, 5 \} \cup \{ p \text{ prime } : G(p) \subsetneq GL_2(\mbz/p\mbz) \; \text{ or } \; p \mid \gD_E \}.
\]
For each prime $p \in S$, define the exponents
\[
\ga_p := \max \left\{ 1, \text{ the exponent $n_\mbq(p)$ of Theorem \ref{araisthm} } \right\}
\]
and
\[
\gb_p := \text{ the exponent of $p$ occurring in } \; \left| GL_2\left(\mbz \Big/ \left( \prod_{\ell \in S \backslash \{ p \} }  \ell \right) \mbz \right) \right|.
\]
Finally, define the positive integer
\begin{equation} \label{defofne}
n_E := \prod_{p \in S} p^{\ga_p + \gb_p}.
\end{equation}
Note that, for $p \in S$ and $M$ dividing $(n_E/p^{\ga_p+\gb_p})^\infty$, one has
\begin{equation} \label{keyfactaboutmp}
\gb_p = \text{ the exponent of $p$ in } |GL_2(\mbz/M\mbz)|.
\end{equation}
Using the above definitions and facts, we will prove
\begin{theorem} \label{nethm}
Let $E$ be any elliptic curve defined over $\mbq$.  Then
\[
\varphi_E(\gal(\ol{\mbq}/\mbq)) = \pi^{-1}(\gal(\mbq(E[n_E])/\mbq)),
\]
where $n_E$ is defined in \eqref{defofne}.  In particular, $m_E \leq n_E$.
\end{theorem}

Note that 
\[
\prod_{p \in S} p^{\gb_p} \leq \left| GL_2 \left(\mbz/ \left( \prod_{\ell \in S }  \ell \right) \mbz \right) \right| \ll \prod_{\ell \in S} \ell^4,
\]
so that, by \eqref{tatemodulesurjection} and \eqref{mazursresult}, if $E$ is semi-stable and non-CM then
\begin{equation} \label{NEbound}
n_E \ll (\prod_{\ell \mid \gD_E} \ell )^5,
\end{equation}
and an affirmative answer to Question \ref{serresquestion} for $K = \mbq$ would imply the above
bound for all non-CM elliptic curves $E$ over $\mbq$.
Thus, Theorem \ref{main} is a corollary of Theorem \ref{nethm}. \\

\noindent \emph{Proof of Theorem \ref{nethm}.}
First we will prove
\begin{lemma} \label{stabilizes}
For any positive integer $n_1$ dividing $n_E^\infty$, one has
\[
G(n_1) \simeq \pi^{-1}(G(d)),
\]
where $d$ is the greatest common divisor of $n_1$ and $n_E$.
\end{lemma}
In the language of \cite{langtrotter}, this lemma says that $n_E$ ``stabilizes'' the Galois representation $\varphi_E$.
The second lemma says that $n_E$ ``splits'' $\varphi_E$ as well.
\begin{lemma} \label{splits}
For any positive integers $n_1$ dividing $n_E^\infty$ and $n_2$ coprime to $n_E$, one has
\[
G(n_1n_2) \simeq G(n_1) \times GL_2(\mbz/n_2\mbz).
\]
\end{lemma}
The two lemmas together imply Theorem \ref{nethm}.
\hfill $\Box$ \\

\noindent \emph{Proof of Lemma \ref{stabilizes}.}
Fix an arbitrary divisor $d$ of $n_E$.  The statement of the lemma is trivial if $n_1 = d$.  Now we will prove it by induction on the set
\[
\mc{N}_d := \{ n \in \mathbb{N} : n \text{ divides } n_E^\infty, \; \gcd(n,n_E) = d \}.
\]
Let $n_1 \in \mc{N}_d$ and suppose that for each $n \in \mc{N}_d \cap \{ 1, 2, \dots, n_1-1\}$, the statement of the lemma is true.  Notice that if $n_1 > d$, then there must exist a prime $p \in S$ satisfying
\[
p^{\ga_p + \gb_p} \text{ exactly divides } d \text{ and } p^{\ga_p + \gb_p + 1} \text{ divides } n_1.
\]
Write $n_1 = p^{r+1} M$, where $p$ does not divide $M$ and 
\begin{equation} \label{r}
r \geq \ga_p+\gb_p.  
\end{equation}
We will show that
\begin{equation} \label{notgrown}
L_{p^{r+1}} \cap L_M = L_{p^r} \cap L_M.
\end{equation}
If this is true, then, writing $k$ for this common field, we have that
\[
\gal(L_{p^{r+1}}L_M/k) \simeq \gal(L_{p^{r+1}}/k) \times \gal(L_M/k) 
\]
and
\[
\gal(L_{p^{r}}L_M/k) \simeq \gal(L_{p^{r}}/k) \times \gal(L_M/k),
\]
from which it follows that $[L_{p^{r+1}M} : L_{p^r}L_M ] = [L_{p^{r+1}} : L_{p^r} ]$.  Since $r \geq \ga_p$, we conclude that
\[
G(n_1) = \pi^{-1}(G(p^r M)),
\]
proving the lemma by induction.

To see why \eqref{notgrown} holds, let us write
\begin{equation} \label{containment}
F_x := L_{p^x} \cap L_M \subseteq L_M \quad \quad \quad (x \geq 1).
\end{equation}
Note that, for $x \geq 1$, the degree $[F_{x+1} : F_{x}]$ is always a power of $p$.  Thus, if $\gb_p = 0$, then by \eqref{keyfactaboutmp}, we must have $F_r = F_{r+1}$.  Now assume that $\gb_p \geq 1$. Suppose first that
\begin{equation*}
\forall s \in \{1, 2, \dots, r - \ga_p \}, \quad F_{\ga_p + s-1} \subsetneq F_{\ga_p + s}.  
\end{equation*}
By \eqref{containment}, \eqref{r}, and \eqref{keyfactaboutmp} we see that this may only happen if $r = \gb_p+\ga_p$ and the exponent of $p$ in $[F_{r} : \mbq]$ is $\gb_p$.
In this case we see from \eqref{containment} that $F_{r+1} = F_r$.  

Now suppose instead that for some $s \in \{ 1, 2, \dots, r - \ga_p \}$ one has $ F_{\ga_p + s-1} = F_{\ga_p + s}$.  We'll first show that under these conditions, $F_{\ga_p + s-1} = F_{\ga_p + s + 1}$.  To ease notation, we will write $\ga := \ga_p + s-1$, so that we are trying to prove that
\[
F_\ga = F_{\ga+1} \Longrightarrow F_\ga = F_{\ga+2}.
\]
Denote by
\[
\pi_2 : G(p^{\ga+2}) \rightarrow G(p^{\ga + 1}), \quad \pi_1 : G(p^{\ga + 1}) \rightarrow G(p^{\ga})
\]
the restrictions of the natural projections and let $N' \subseteq N \subseteq G(p^{\ga+2})$ be the normal subgroups satisfying
\[
F_{\ga} = F_{\ga+1} = L_{p^{\ga+2}}^N \quad \text{and} \quad F_{\ga+2} = L_{p^{\ga+2}}^{N'}.
\]
Our contention is that $N' = N$.  Now,
\begin{equation} \label{kerinN}
L_{p^{\ga+2}}^{ \ker \pi_2 \cdot N' } = L_{p^{\ga+2}}^{\ker \pi_2} \cap L_{p^{\ga+2}}^{N'} = L_{p^{\ga+2}}^N,
\end{equation}
which implies that the restriction of $\pi_2$ to $N'$ maps surjectively onto $\pi_2(N)$:
\[
N' \twoheadrightarrow \pi_2(N).
\]
The fact that $L_{p^{\ga+2}}^N = F_\ga \subseteq L_{p^\ga} = L_{p^{\ga+2}}^{\ker (\pi_1 \circ \pi_2)}$ implies that
\[
\pi_2^{-1}(\ker \pi_1) = \ker (\pi_1 \circ \pi_2) \subseteq N \subseteq \pi_2^{-1} (\pi_2(N)),
\]
so that 
\begin{equation*} \label{kerpicontainedinpi}
\ker \pi_1 \subseteq \pi_2(N).
\end{equation*}
Since $\ga \geq \ga_p$, we know that
\[
\ker \pi_2 = I + p^{\ga+1} M_{2 \times 2}(\mbz/p\mbz) \quad \text{and} \quad \ker \pi_1 = I + p^{\ga}M_{2 \times 2}(\mbz/p\mbz). 
\]
Now pick any
\[
I + p^\ga A \in \ker \pi_1
\]
and find a pre-image $X = I + p^\ga A + p^{\ga+1}B \in N'$.  But then
\[
X^p \equiv I + p^{\ga+1}A \mod p^{\ga+2} \in N',
\]
and so $I + p^{\ga+1}M_{2\times 2}(\mbz/p\mbz) = \ker \pi_2 \subseteq N'$.  This together with \eqref{kerinN} shows that $N' = N$, as desired.  Replacing $s$ by $s+1$ and repeating the argument inductively, we conclude that $F_{\ga_p+s-1} = F_{\ga_p + k}$ for any positive integer $k \geq s-1$, so that in particular $F_{r+1} = F_r$.  This finishes the proof of Lemma \ref{stabilizes}.
\hfill $\Box$ \\

\noindent \emph{Proof of Lemma \ref{splits}.}
The reasoning here is very similar to that of \cite[Theorem $6.1$, p. $49$]{langtrotter}. The first step is to prove
\begin{sublemma} \label{sublemma}
Fix any integers $M_1$ and $M_2$ with the property that $2 \nmid M_2$, $5 \nmid M_2$, and $\gcd(M_1 \gD_E ,M_2) = 1$.  If $G(M_2) = GL_2(\mbz/M_2\mbz)$, then
\[
G(M_1M_2) \simeq G(M_1) \times GL_2(\mbz/M_2\mbz).
\]
\end{sublemma}
\noindent \emph{Proof of Sublemma \ref{sublemma}.}  
Set $F := L_{M_1} \cap L_{M_2}$.  We need to show that $F = \mbq$.  Suppose that $F \neq \mbq$.  Note that $1 \neq \gal(F/\mbq)$ is a common quotient group of $G(M_1)$ and $G(M_2) = GL_2(\mbz/M_2\mbz)$.  Replacing $F$ by a subfield, we may assume that $\gal(F/\mbq)$ is a common non-trivial \emph{simple} quotient.  We claim that this common simple quotient must be abelian.  For a finite group $G$ let $\text{Occ}(G)$ denote the set of simple non-abelian groups which occur as quotients of subgroups of $G$.  One easily deduces from \cite[p. IV-25]{serre2} that, for any positive integer $M$, $\text{Occ}(GL_2(\mbz/M\mbz))$ is equal to
\[
\left( \bigcup_{{\begin{substack} {p \mid M \\ p > 5 \\ p \equiv \pm 1 \; \text{mod } 5} \end{substack}}} \{ PSL_2(\mbz/p\mbz) , A_5 \} \right) \cup \left( \bigcup_{{\begin{substack} {p \mid M \\ p > 5 \\ p \equiv \pm 2 \; \text{mod } 5} \end{substack}}} \{ PSL_2(\mbz/p\mbz) \} \right) \cup \left( \bigcup_{{\begin{substack} {p \mid M \\ p = 5 } \end{substack}}} \{ A_5 \} \right).
\]
(Note that $A_5 \simeq PSL_2(\mbz/5\mbz)$.)  One can use elementary group theory to show that 
\[
\left\{ \text{ simple non-abelian quotients of } GL_2(\mbz/M\mbz) \right\} \subseteq \bigcup_{{\begin{substack} {p \mid M \\ p > 3 } \end{substack}}} \{ PSL_2(\mbz/p\mbz) \}.
\]
Thus, the assumptions on $M_1$ and $M_2$ imply that $\gal(F/\mbq)$ must be abelian.  Since $M_2$ is odd, the commutator subgroup 
\[
[GL_2(\mbz/M_2\mbz),GL_2(\mbz/M_2\mbz)] = SL_2(\mbz/M_2\mbz), 
\]
which implies that $F$ is contained in the cyclotomic field
\[
F \subseteq \mbq \left(\exp \left( \frac{2\pi i}{M_2} \right) \right).
\]
Let $p$ be a prime ramified in $F$.  We see that $p$ must divide the discriminants of both $L_{M_1}$ and $\mbq \left(\exp \left( \frac{2\pi i}{M_2} \right) \right)$, which is impossible since $\gcd(M_1 \gD_E ,M_2) = 1$.  Since $\mbq$ has no everywhere unramified extensions, we have arrived at a contradiction.  Thus, we cannot have $F \neq \mbq$, and the sublemma is proved.
\hfill $\Box$ \\

To prove Lemma \ref{splits}, we first prove by induction on the number of primes $p$ dividing $n_2$, that in fact
\begin{equation} \label{fullatN2}
G(n_2) \simeq GL_2(\mbz/n_2\mbz).
\end{equation}
The case where $n_2$ is a power of a prime $p > 5$ follows from \eqref{tatemodulesurjection}.
Then, \eqref{fullatN2} is proved by writing $n_2 = p^n M$ with $n \geq 1$ and $p \nmid M$ and applying Sublemma \ref{sublemma} with $M_1 = p^n$ and $M_2 = M$.  Finally, to prove Lemma \ref{splits}, we apply the sublemma with $M_i = n_i$.
\hfill $\Box$

We end by asking the following weakening of Question \ref{serresquestion}.
\begin{question}
Fix a number field $K$.  Does there exist a constant $C_K$ so that for each prime number $p$ one has
\[
n_K(p) \leq C_K,
\]
where $n_K(p)$ is the exponent occurring in Theorem \ref{araisthm}?
\end{question}
Conditional upon an affirmative answer to this question, Theorem \ref{nethm} together with \cite[Theorem $2$]{cojocarukani} would imply that, for any non-CM elliptic curve $E$ over $\mbq$ one has
\[
m_E \ll \left( \prod_{p \leq B_E} p \right)^{C_\mbq + 4} \cdot \left( \prod_{p \mid \gD_E} p \right)^5,
\]
where 
\[
B_E := \frac{4\sqrt{6}}{3} \cdot N_E \prod_{p \mid \gD_E} \left( 1 + \frac{1}{p} \right)^{1/2} + 1,
\]
$N_E$ denoting the conductor of $E$.

\vspace{.5in}

\begin{center}
Centre de Recherches Math\'{e}matiques \\
Universit\'{e} de Montr\'{e}al \\
P.O. Box 6128, \\
Centre-ville Station \\
Montr\'{e}al, Qu\'{e}bec  H3C 3J7, Canada. \\
E-mail:  jones@dms.umontreal.ca
\end{center}

\end{document}